\documentclass[runningheads]{llncs}

\usepackage[T1]{fontenc}
\usepackage{amsmath,amssymb,mathtools}
\usepackage{booktabs,array}
\usepackage{graphicx}
\usepackage{tikz}
\usetikzlibrary{arrows.meta,positioning,calc,shapes.geometric}
\usepackage{microtype}

\newcommand{\R}{\mathbb{R}}
\newcommand{\C}{\mathbb{C}}
\newcommand{\Z}{\mathbb{Z}}
\newcommand{\Oa}{\mathbb{O}}
\newcommand{\pO}{p\mathbb{O}}
\newcommand{\Ok}{\mathcal{O}}
\newcommand{\OEight}{\mathbb{O}_{E_8}}

\newcommand{\e}{\mathfrak{e}}
\newcommand{\slthree}{\mathfrak{sl}_3}
\newcommand{\Idem}{\operatorname{Idem}}
\newcommand{\Aut}{\operatorname{Aut}}
\newcommand{\Fix}{\operatorname{Fix}}

\newcommand{\Ree}{\operatorname{Re}}
\newcommand{\Span}{\operatorname{Span}}
\newcommand{\PhiE}{\Phi(E_8)}
\newcommand{\ip}[2]{\langle #1,#2\rangle}
\definecolor{istblue}{HTML}{009DE0}
\definecolor{istbluemd}{HTML}{0079A8}
\definecolor{istbluedk}{HTML}{004965}
\definecolor{istgrey}{HTML}{46555F}
\definecolor{istgreylt}{HTML}{D2D3D4}
\definecolor{istink}{HTML}{3A454C}
\colorlet{okubonavy}{istbluedk}   
\colorlet{okubosky}{istbluemd}    
\colorlet{okubogold}{istblue}     
\colorlet{okubogolddk}{istbluemd} 
\colorlet{okuboink}{istink}       
\colorlet{okuboskyL}{istblue!60}  
\begin{document}

\title{Integral Orders for the Okubo Algebra and Idempotent Geometries of the $E_8$ Lattice}
\titlerunning{Integral Okubo Orders and Idempotent Geometries of $E_8$}

\author{
Daniele Corradetti\inst{1}
\and
Alessio Marrani\inst{2}
}

\authorrunning{D. Corradetti and A. Marrani}

\institute{
Grupo de F\'isica Matem\'atica, Instituto Superior T\'ecnico,
Universidade de Lisboa,\\
Av. Rovisco Pais, 1049-001 Lisboa, Portugal\\
\email{danielecorradetti@tecnico.ulisboa.pt}
\and
Centre for Mathematics and Theoretical Physics,
University of Hertfordshire,\\
Hatfield AL10 9AB, UK\\
\email{a.marrani@herts.ac.uk}
}

\maketitle

\begin{abstract}
We study the Coxeter--Dickson $E_8$ order as a common integral support for the octonionic, para-octonionic and compact real Okubo products.  The three algebra structures have the same additive lattice, the same positive composition norm and hence the same $240$ norm-one vectors, namely the roots of $E_8$ and the vertices of the Gosset polytope $4_{21}$.  Their multiplicative structures are nevertheless different and this difference is already visible in the idempotents: among the same $240$ roots one finds respectively $1$, $57$ and $12$ nonzero integral idempotents.  We show that these three counts organize the common root system in three increasingly polarized ways.  The octonionic product leaves the full $E_8$ geometry unseparated; the $57$ para-octonionic idempotents reproduce the $E_7$ contact decomposition $1+56+126+56+1$; the $12$ Okubo idempotents split into four mutually orthogonal oriented $A_2$ triangles and hence determine canonically an $A_2^4$ subsystem.  Choosing one triangle as the external $A_2$ then yields the $A_2+E_6$ Magic Star, while the remaining nine idempotents determine the trinification subsystem $A_2^3\subset E_6$.  This provides an arithmetic interpretation of the $E_8\supset E_7$ and $E_8\supset E_6+A_2$ decompositions directly from integral idempotents.

\keywords{Okubo algebra \and octonions \and para-octonions \and $E_8$ lattice \and idempotents \and Magic Star \and Jordan algebra \and exceptional Lie algebras}
\end{abstract}

\section{Introduction and motivations}

Integral numbers in division algebras lie at a classical intersection of algebra, arithmetic and crystallographic geometry \cite{Corradetti2026}.  The Gaussian and Eisenstein integers, the Hurwitz quaternions and the Coxeter--Dickson octonions provide progressively richer examples in dimensions $2$, $4$ and $8$ \cite{Dickson,Coxeter,ConwaySmith,Johnson}. The insights of N. Johnson elucidated two main ingredients that a ring needs in order to be a basic system of integers, i.e. to be in a lattice and to be embedded in a composition algebra. Perhaps the most striking result on integral systems over composition algebras is given by the octonionic case. Indeed, it is very well known that the maximal order $\OEight$ is isometric to the exceptional lattice $E_8$; its $240$ norm-one elements are the $240$ roots of $E_8$ and, equivalently, the vertices of the Gosset polytope $4_{21}$. The eight-dimensional real composition division algebras, however, are not exhausted by the ordinary octonions \cite{Corr24-3DGA,Corr26-8DComposition}.  Besides the unital octonions $\Oa$ one has the para-octonions $\pO$, which are para-unital, and the compact real Okubo algebra $\Ok$, which has neither a unit nor a paraunit and whose automorphism group is $PSU(3)$ rather than $G_2$ \cite{Okubo,ElduqueMyung,Elduque2015,CZ22-OkuboSpin,MCZ25-OkuboPhysics}.  These three algebras are related by isotopy, but they are not isomorphic as algebras.

The arithmetic observation underlying this note is that the same Coxeter--Dickson lattice can support all three products.  For the para-octonions this is immediate from closure under conjugation and octonionic multiplication.  For the Okubo case it is less obvious since one has to carefully choose the automorphism $\tau\in G_2$ for the Petersson product \cite{Corradetti2026} in order to have a realisation of the compact real Okubo algebra that also preserves the four crystallographic octonionic orders with metric lattices
\[
C_8,\qquad A_2^4,\qquad D_4^2,\qquad E_8.
\]

 Once the metric lattice is fixed, the multiplication does not move the $240$ roots; it only changes how they multiply.  The idempotent equation therefore becomes a particularly economical multiplicative invariant.  On $E_8$ the three products give
\[
\boxed{1,\qquad 57,\qquad 12}
\]
integral idempotents.  We shall see that these three numbers are not merely counts. In fact, they reveal three different exceptional organizations of the same root system, centered respectively on $E_8$, $E_7$ and $E_6$.

\section{From integral elements to the $E_8$ shell}

What is an integer? This apparently simple question has a not-so-simple answer. Following Dickson, given
a composition algebra $\mathbb{A}$, i.e. an algebra with a norm $n$
such that $n\left(x\cdot y\right)=n\left(x\right)n\left(y\right)$,
endowed with an involution $x\longmapsto\overline{x}$ called
\emph{conjugation}, then a set $I$ of \emph{integral elements} of
$\mathbb{A}$ is a set that is closed under addition and multiplication
that contains the element $1$ and such that for every element $x\in I$
we have 
\begin{align}
\operatorname{tr}\left(x\right)=x+\overline{x} & \in\mathbb{Z},\label{eq:trace1}\\
n\left(x\right)=x\overline{x} & \in\mathbb{Z}.\label{eq:norm1}
\end{align}
Although this definition is broad and was widely accepted, it has two
notable shortcomings. First, it can easily yield an infinite number
of integral elements over Hurwitz division algebras.
For example, for any natural number $m\in\mathbb{N}$, the subring
$\mathbb{Z}\left[\sqrt{-m}\right]$ is a set of integral elements
of $\mathbb{C}$. Second, its specificity, requiring the presence
of the unit element, automatically excludes integral systems from
non-unital algebras, such as para-octonions $p\mathbb{O}$ and Okubo
algebra $\mathcal{O}$.

The first issue with Dickson's definition was addressed by Johnson
\cite{Johnson} who had the intuition that integers, at the end of the day, are lattices embedded in composition algebras. He then
arrived at the following definition: A \emph{basic system of integral
elements} over a real algebra is given by a set of elements such that:
\begin{enumerate}
\item The trace and the norm are integers. Specifically,
\begin{align*}
\operatorname{tr}\left(x\right) & =x+\overline{x}\in\mathbb{Z},\\
n\left(x\right) & =x\overline{x}\in\mathbb{Z},
\end{align*}
 and where the conjugation is given by an involution of the algebra.
\item The elements form a subring of the algebra, closed under multiplication,
and comprising a set of invertible unit elements.
\item The elements span a two-, four-, or eight-dimensional lattice embedded
in $\mathbb{R},\mathbb{C},\mathbb{H}$ or $\mathbb{O}$.
\end{enumerate}

Both formulations presuppose a unit, so neither applies verbatim to the para-octonions or to the Okubo algebra. For the non-unital products of Section~3 we therefore use the following neutral notion: given an eight-dimensional real composition algebra $\left(A,\circ,n\right)$ and a nonzero idempotent $e\in A$, a \emph{$\mathbb{Z}$-integral system relative to $e$} is a free $\mathbb{Z}$-module $L\subset A$ of rank eight with $e\in L$, $L\circ L\subseteq L$, $n\left(L\right)\subseteq\mathbb{Z}$ and $\ip{L}{e}\subseteq\mathbb{Z}$. For $A$ unital and $e=1$ this reduces to the definition above, and we reserve the word \emph{order} for that case.

Taking into account Johnson's refinements, the number of viable sets
of integral elements is reduced to ten (1 for the reals $\mathbb{R}$, 2 for the complex numbers $\mathbb{C}$, 3 for the quaternions $\mathbb{H}$ and 4 for the octonions $\mathbb{O}$), as detailed in Table
\ref{tab:Johnson Integers}.
\begin{table}
\centering{}%
\begin{tabular}{|c|c|c|c|c|c|}
\hline 
\textbf{Name} & \textbf{Alg.} & \textbf{Dim.} & \textbf{Symbol} & \textbf{Unit El.} & \textbf{Lattice}\tabularnewline
\hline 
\hline 
Integers & $\mathbb{R}$ & 1 & $\mathbb{Z}$ & 2 & $A_{1}$\tabularnewline
\hline 
Eisenstein & $\mathbb{C}$ & 2 & $\mathbb{C}_{A_{2}}$ & 6 & $A_{2}$\tabularnewline
\hline 
Gaussian & $\mathbb{C}$ & 2 & $\mathbb{C}_{C_{2}}$ & 4 & $C_{2}$\tabularnewline
\hline 
Hamilton & $\mathbb{H}$ & 4 & $\mathbb{H}_{2C_{2}}$ & 8 & $C_{2}\oplus C_{2}$\tabularnewline
\hline 
Hybrid & $\mathbb{H}$ & 4 & $\mathbb{H}_{2A_{2}}$ & 12 & $A_{2}\oplus A_{2}$\tabularnewline
\hline 
Hurwitz & $\mathbb{H}$ & 4 & $\mathbb{H}_{D_{4}}$ & 24 & $D_{4}$\tabularnewline
\hline 
Cayley-Graves & $\mathbb{O}$ & 8 & $\mathbb{O}_{C_{8}}$ & 16 & $C_{8}$\tabularnewline
\hline 
Comp. Eisenstein & $\mathbb{O}$ & 8 & $\mathbb{O}_{4A_{2}}$ & 24 & $A_{2}\oplus A_{2}\oplus A_{2}\oplus A_{2}$\tabularnewline
\hline 
Coupled Hurwitz & $\mathbb{O}$ & 8 & $\mathbb{O}_{2D_{4}}$ & 48 & $D_{4}\oplus D_{4}$\tabularnewline
\hline 
Coxeter-Dickson & $\mathbb{O}$ & 8 & $\mathbb{O}_{E_{8}}$ & 240 & $E_{8}$\tabularnewline
\hline 
\end{tabular}\caption{\label{tab:Johnson Integers}Summary of all sets
of integral elements over the division Hurwitz algebras $\mathbb{R},\mathbb{C},\mathbb{H}$
and $\mathbb{O}$. In the first column we indicated the name according
to \cite{Johnson}; then the related algebra in which the integral set
is embedded; the dimension of the algebra; a notational symbol we
introduced; the number of invertible unit elements; finally,
in the last column the lattice associated with the integral set.}
\end{table}

The octonionic case is the one of interest here. The \emph{Cayley-Graves integers} $\mathbb{O}_{C_{8}}$ are the simplest
octonionic integral set, consisting of octonions with all-integer coefficients
in the standard basis $\left\{ 1,\mathrm{i},\mathrm{j},\mathrm{k},\mathrm{l},\mathrm{il},\mathrm{jl},\mathrm{kl}\right\} $.
They possess 16 unit elements forming the Moufang loop $M_{16}(Q_{8})$,
and their lattice is $C_{8}$. The \emph{Compounded Eisenstein integers} $\mathbb{O}_{4A_{2}}$ are
obtained by embedding four copies of the Eisenstein integers into $\mathbb{O}$
using compatible quaternionic subalgebras. They possess 24 unit elements
forming the loop $M_{24}(4A_{2})$, and their lattice is $A_{2}\oplus A_{2}\oplus A_{2}\oplus A_{2}$. The \emph{Coupled Hurwitz integers} $\mathbb{O}_{2D_{4}}$ are obtained
by coupling two copies of the Hurwitz quaternionic integers in a compatible
manner. They possess 48 unit elements forming the loop $M_{48}(2D_{4})$,
and their lattice is $D_{4}\oplus D_{4}$.
Finally, the maximal integral set is the \emph{Coxeter-Dickson octonions}
$\mathbb{O}_{E_{8}}$, which are the elements of the form

\begin{equation}
x=a_{0}+a_{1}\text{i}+a_{2}\text{j}+a_{3}\text{k}+a_{4}\text{h}+a_{5}\text{ih}+a_{6}\text{jh}+a_{7}\text{kh},
\end{equation}
 where $\text{h}=\left(\text{i}+\text{j}+\text{k}+\text{l}\right)/2$
and $a_{0},...,a_{7}\in\mathbb{Z}$. 

Indeed, it is easy to see that these elements are closed under octonionic conjugation and multiplication and that the norm and the trace are integers. Throughout, the octonions are endowed with the Euclidean bilinear form polarizing the composition norm,
\begin{equation}
\ip{x}{y}:=x\overline{y}+y\overline{x}=2\Ree\left(x\overline{y}\right),\qquad\ip{x}{x}=2n\left(x\right),\label{eq:normalisation}
\end{equation}
so that elements of composition norm one have squared length $2$. With this normalisation, the isometry between the set of integral octonions $\mathbb{O}_{E_{8}}$
and the $E_{8}$-lattice is a
well-known result following from a straightforward calculation. We denote by $\OEight(1)$ the shell of elements of unit norm: these are $240$ in number, they form the root system $\PhiE$ and they are the vertices of the Gosset polytope $4_{21}$, whose symmetry group is of course the Weyl group of $E_8$.
The unit invertible elements of the set are the 240 elements 

\begin{equation}
\begin{array}{cc}
\pm1,\pm\text{i},\pm\text{j},\pm\text{k}, & \pm\text{l},\pm\text{i\text{l}},\pm\text{j\text{l}},\pm\text{k\text{l}},\\
\frac{1}{2}\left(\pm1\pm\text{j}\pm\text{k}\pm\text{i\text{l}}\right), & \frac{1}{2}\left(\pm\text{i}\pm\text{l}\pm\text{j\text{l}}\pm\text{k\text{l}}\right),\\
\frac{1}{2}\left(\pm1\pm\text{k}\pm\text{i}\pm\text{j\text{l}}\right), & \frac{1}{2}\left(\pm\text{j}\pm\text{l}\pm\text{k\text{l}}\pm\text{i\text{l}}\right),\\
\frac{1}{2}\left(\pm1\pm\text{i}\pm\text{j}\pm\text{k\text{l}}\right), & \frac{1}{2}\left(\pm\text{k}\pm\text{l}\pm\text{i\text{l}}\pm\text{j\text{l}}\right),\\
\frac{1}{2}\left(\pm1\pm\text{i\text{l}}\pm\text{j\text{l}}\pm\text{k\text{l}}\right), & \frac{1}{2}\left(\pm\text{i}\pm\text{j}\pm\text{k}\pm\text{l}\right),\\
\frac{1}{2}\left(\pm1\pm\text{i}\pm\text{l}\pm\text{i\text{l}}\right), & \frac{1}{2}\left(\pm\text{j}\pm\text{k}\pm\text{j\text{l}}\pm\text{k\text{l}}\right),\\
\frac{1}{2}\left(\pm1\pm\text{j}\pm\text{l}\pm\text{j\text{l}}\right), & \frac{1}{2}\left(\pm\text{k}\pm\text{i}\pm\text{k\text{l}}\pm\text{i\text{l}}\right),\\
\frac{1}{2}\left(\pm1\pm\text{k}\pm\text{l}\pm\text{k\text{l}}\right), & \frac{1}{2}\left(\pm\text{i}\pm\text{j}\pm\text{i\text{l}}\pm\text{j\text{l}}\right),
\end{array}\label{eq:InvertibleElements}
\end{equation}
In fact, the elements in (\ref{eq:InvertibleElements}) form a well-known
Moufang loop.

\section{Octonions, para-octonions and the integral Okubo product}

Let $(\Oa,\cdot,n)$ denote the usual algebra of octonions. Then one can define the para-octonionic product as
\begin{equation}
 x\bullet y=\bar x\cdot\bar y .
 \label{eq:para}
\end{equation}
Since the octonions are a composition algebra and $n(\bar x)=n(x)$, one has that also this new algebra is a composition algebra, i.e., $n(x\bullet y)=n(x)n(y)$, but instead of a unit it has a paraunit $1$, in the sense that
\[
1\bullet x=x\bullet1=\bar x.
\]
In particular $1$ is an idempotent since $1\bullet1=1$, although $1$ is not an identity for $\bullet$.

For the Okubo product we fix the octonionic multiplication on the basis $\left\{ 1,\text{i},\text{j},\text{k},\text{l},\text{il},\text{jl},\text{kl}\right\}$ by the seven oriented Fano lines
\begin{equation}
\begin{array}{c}
\left(\text{i},\text{j},\text{k}\right),\quad\left(\text{i},\text{l},\text{il}\right),\quad\left(\text{j},\text{l},\text{jl}\right),\quad\left(\text{k},\text{l},\text{kl}\right),\\[2pt]
\left(\text{jl},\text{i},\text{kl}\right),\quad\left(\text{kl},\text{j},\text{il}\right),\quad\left(\text{il},\text{k},\text{jl}\right),
\end{array}\label{eq:fano}
\end{equation}
where $\left(x,y,z\right)$ stands for $xy=z$, $yz=x$, $zx=y$, and we relabel the imaginary units by
\begin{equation}
\left(e_{1},\ldots,e_{7}\right)=\left(\text{i},\text{j},\text{l},\text{k},\text{kl},\text{jl},\text{il}\right).\label{eq:fanobasis}
\end{equation}
Consider then the order-three signed permutation
\begin{align}
&\tau(1)=1,\qquad \tau(e_3)=e_3,\nonumber\\
&\tau(e_1)=e_2,\quad \tau(e_2)=-e_4,\quad \tau(e_4)=-e_1,\nonumber\\
&\tau(e_5)=-e_7,\quad \tau(e_7)=e_6,\quad \tau(e_6)=-e_5.
\label{eq:tau}
\end{align}
Inspection of the oriented Fano lines \eqref{eq:fano} gives $\tau(xy)=\tau(x)\tau(y)$, while the two signed three-cycles give $\tau^3=\mathrm{id}$.  Hence $\tau\in G_2$ is an integral octonion automorphism of order three.  Its fixed subalgebra is quaternionic: it is spanned by $1$, $e_3$ and the two $\tau$-invariant vectors
\begin{equation}
a=-e_{1}-e_{2}+e_{4},\qquad b=-e_{5}+e_{6}+e_{7},\label{eq:fixedquat}
\end{equation}
which satisfy $a^{2}=b^{2}=-3$ and $ab=3e_{3}=-ba$, so that $\Fix(\tau)=\Span_{\R}\left\{ 1,e_{3},a,b\right\}$. Thus, by the classification of order-three Petersson algebras \cite{Petersson,Elduque2018}, the isotope
\begin{equation}
 x*y=\tau(\bar x)\cdot\tau^2(\bar y)
 \label{eq:okubo}
\end{equation}
is the compact real Okubo division algebra \cite{Elduque2018}.

The crucial arithmetic point is that $\tau$ preserves the four octonionic orders.  Since each order is also closed under conjugation and octonionic multiplication, it is closed under both \eqref{eq:para} and \eqref{eq:okubo}.  The result is that one obtains three integral algebra structures on the same additive $E_8$ lattice.
\begin{proposition}[One lattice, three products]
The identity map on $\OEight$ is an isometry among the octonionic, para-octonionic and Okubo integral systems.  In particular, all three have the same Gram matrix, the same $240$ norm-one vectors and the same Gosset polytope $4_{21}$.  They are pairwise non-isomorphic as algebras.
\end{proposition}
\noindent
The difference appears only after multiplication is introduced.  With the ordinary product the roots form the Moufang loop $M_{240}(E_8)$; with the para-product they form a para-isotopic quasigroup; with the Okubo product they form an Okubo quasigroup.  The most economical distinction is the idempotent equation.

\section{Three idempotent geometries of the same $E_8$}

Before comparing the three products, note that for any of them a nonzero idempotent $e$ satisfies $n(e)=n(e)^2$ by composition, whence $n(e)=1$ by positive definiteness: every nonzero integral idempotent lies on the common shell $\OEight(1)$, so that testing the $240$ roots is exhaustive.

All root-system decompositions and representation branchings below are written for the complexified exceptional Lie algebras, whereas statements about $G_2$, $PSU(3)$ or $SU(3)^3$ refer to the corresponding compact real forms; the root combinatorics does not depend on this choice.

In a division algebra an idempotent $e$ satisfies $e^2=e$ and, if such an algebra is also unital, then the only nonzero solution is
\[
e=1.
\]
This is the case for the octonionic algebra, which is unital and division. Thus the octonionic product singles out the multiplicative identity but does not select a nontrivial oriented subsystem among the $240$ roots. Things are different for the para-octonionic and the Okubo cases.

Let us consider the \emph{para-octonionic case}.
For a norm-one element $x\in\OEight(1)$ the para-idempotent equation $x\bullet x=x$ reads $\overline{x}^{\,2}=x$, equivalently $x^{2}=\overline{x}$.
Since the subalgebra generated by one octonion is associative and $\bar x=x^{-1}$, this is equivalent to $x^3=1$.  Hence either $x=1$, or $x^2+x+1=0,$ which for a unit octonion is equivalent to
\[
\Ree(x)=-\frac12.
\]
There are exactly $56$ Coxeter--Dickson roots with this real part.  Therefore since
\begin{equation}
\Idem_{p\Oa}(\OEight(1))=\{1\}\sqcup\{x\in\PhiE:\ip{x}{1}=-1\},
\label{eq:57}
\end{equation}
one has that the total number of idempotents on the lattice for the para-octonionic case is $|\Idem_{p\Oa}|=57$.

Now regard the para-unit $\alpha=1$ as a root.  Relative to $\alpha$, the $240$ roots are distributed according to the possible inner products $2,1,0,-1,-2$ as
\begin{equation}
\boxed{240=1+56+126+56+1.}
\label{eq:e7roots}
\end{equation}
This is the standard root-layer distribution about a root of $E_8$ \cite{Adams}.
Indeed, the $126$ roots orthogonal to $\alpha$ form a root system of type $E_7$. More strikingly, the corresponding Lie-algebra decomposition is the well-known contact $5$-grading
\begin{equation}
\e_8=\mathbf1_{-2}\oplus\mathbf{56}_{-1}\oplus(\e_7\oplus\C)_0
\oplus\mathbf{56}_{+1}\oplus\mathbf1_{+2}.
\label{eq:contact}
\end{equation}
At the root level one of the two outer halves of \eqref{eq:e7roots} is precisely the para-idempotent set: the distinguished point $1$ together with the $56$ roots at inner product $-1$.

Reducing the integer degrees in \eqref{eq:contact} modulo three gives a cyclic grading
\[
\e_8=\mathcal G_0\oplus\mathcal G_1\oplus\mathcal G_2,
\]
with
\[
\mathcal G_0=\e_7\oplus\C,
\quad
\mathcal G_1=\mathbf{56}_{+1}\oplus\mathbf1_{-2},
\quad
\mathcal G_2=\mathbf{56}_{-1}\oplus\mathbf1_{+2}.
\]
Thus the number $57$ has a direct grading interpretation: it is the root count in either nonzero cyclic component.

The $56$ nontrivial para-idempotents are the weights occurring in the degree $\pm1$ components of the contact grading; the associated root spaces form the minuscule $\mathbf{56}$ of $E_7$.  This representation is the Freudenthal module associated with the Albert algebra and carries an invariant symplectic form together with the exceptional quartic invariant \cite{Slansky}.  Under the subgroup $E_6\times U(1)$ it further decomposes as
\begin{equation}
\mathbf{56}=\mathbf1_{-3}\oplus\mathbf{27}_{-1}\oplus\overline{\mathbf{27}}_{+1}\oplus\mathbf1_{+3},
\label{eq:56e6}
\end{equation}
while
\begin{equation}
\e_7=\overline{\mathbf{27}}_{-2}\oplus(\e_6\oplus\C)_0\oplus\mathbf{27}_{+2}.
\label{eq:e7e6}
\end{equation}
Equations \eqref{eq:56e6}--\eqref{eq:e7e6} provide the bridge from the para-octonionic $E_7$ view to the Jordan-theoretic $E_6$ view appearing below in the Okubo case.

\begin{table}[t]
\centering
\caption{The same $E_8$ shell seen through three products.  The column ``exceptional core'' describes the root subsystem exposed by the idempotent pattern.}
\label{tab:threeviews}
\setlength{\tabcolsep}{3pt}\small\begin{tabular}{p{2.25cm}p{2.05cm}p{1.55cm}p{2.55cm}p{2.15cm}}
\toprule
Product & Unit structure & Idem. & Root organization & Core\\
\midrule
Octonionic & unit $1$ & $1$ & full $240$ & $E_8$\\
Para-octonionic & paraunit $1$ & $57$ & $1+56+126+56+1$ & $E_7$\\
Okubo & neither & $12$ &$A^4_2$ then $6+72+6\cdot27$ & $E_6$\\
\bottomrule
\end{tabular}
\end{table}

\begin{figure}[t]
\centering
\begin{tikzpicture}[scale=0.82, every node/.style={font=\bfseries}]

  \coordinate (T)  at (90:3.05);
  \coordinate (BR) at (-30:3.05);
  \coordinate (BL) at (210:3.05);
  \coordinate (B)  at (270:3.05);
  \coordinate (TR) at (30:3.05);
  \coordinate (TL) at (150:3.05);

  \draw[gray!55,line width=1.05pt] (T)--(BR)--(BL)--cycle;
  \draw[gray!55,line width=1.05pt] (B)--(TR)--(TL)--cycle;

  \foreach \p in {T,TR,BR,B,BL,TL}{
    \node[circle,fill=black,text=white,minimum size=0.56cm,inner sep=0pt] at (\p) {1};
  }

  \foreach \ang in {90,30,-30,-90,-150,150}{
    \node[circle,draw=okubonavy!75,fill=okuboskyL!30,text=okubonavy,
          minimum size=0.72cm,inner sep=0pt] at (\ang:1.62) {27};
  }

  \node[circle,draw=okubogolddk!75,fill=okubogold!18,text=okubonavy,
        minimum size=1.34cm,inner sep=0pt,align=center] at (0,0)
        {$E_6$\\[-0.15em]$72$};


\end{tikzpicture}
\caption{The $A_2+E_6$ Magic Star obtained from an Okubo idempotent triangle: the six external roots of the chosen $A_2$, the $72$ roots of the central $E_6$ and the six $27$-root fibers.}
\label{fig:rootpictures}
\end{figure}
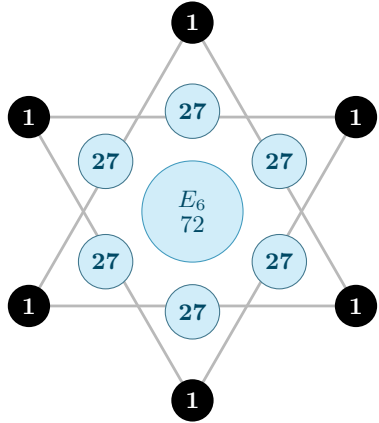

Let us now consider the Petersson product \eqref{eq:okubo}. An integral idempotent is a root $e\in\OEight(1)$ satisfying
\[
\tau(\bar e)\tau^2(\bar e)=e.
\]
\begin{proposition}[Twelve integral Okubo idempotents]
\label{prop:twelve}
With the conventions \eqref{eq:fanobasis} and \eqref{eq:tau}, the exact test of all $240$ roots of $\OEight(1)$ returns exactly twelve solutions, listed in Table~\ref{tab:twelve}.  They split as
\begin{equation}
I_{\rm Ok}=T_1\sqcup T_2\sqcup T_3\sqcup T_4,
\qquad |T_i|=3,
\label{eq:fourtriangles}
\end{equation}
where for every $i$
\[
\sum_{e\in T_i}e=0,
\qquad
\ip{e}{e'}=-1\quad(e\neq e'\in T_i),
\qquad
\ip{e}{f}=0\quad(e\in T_i,\ f\in T_j,\ i\neq j).
\]
\end{proposition}
\noindent
The verification is a finite exact computation over $\Z[\tfrac12]$ in the basis \eqref{eq:fanobasis} with the multiplication \eqref{eq:fano}.  Each $T_i$ is one of the two orientations of a root system
\[
T_i\cup(-T_i)=\Phi(A_2^{(i)}),
\qquad
{i=1,2,3,4},
\]
and, the four being mutually orthogonal, the twelve idempotents canonically determine a rank-eight subsystem $A_2^4\subset\PhiE$.  The first triangle lies in $\Fix(\tau)$, its two nontrivial elements being $\tfrac12\left(-1\pm b\right)$ with $b$ as in \eqref{eq:fixedquat}, and one reads off from Table~\ref{tab:twelve}
\begin{equation}
\tau|_{T_1}=\mathrm{id},
\qquad
T_2\xrightarrow{\tau}T_3\xrightarrow{\tau}T_4\xrightarrow{\tau}T_2.
\label{eq:taucircles}
\end{equation}

\begin{table}[t]
\centering
\caption{The twelve integral Okubo idempotents in the basis \eqref{eq:fanobasis}.  Each row is an oriented $A_2$ triangle; distinct rows are mutually orthogonal.}
\label{tab:twelve}
\setlength{\tabcolsep}{5pt}\small
\begin{tabular}{cccc}
\toprule
$T_1$ & $1$ & $\tfrac12(-1-e_5+e_6+e_7)$ & $\tfrac12(-1+e_5-e_6-e_7)$\\
$T_2$ & $-e_1$ & $\tfrac12(e_1+e_3+e_5+e_6)$ & $\tfrac12(e_1-e_3-e_5-e_6)$\\
$T_3$ & $-e_2$ & $\tfrac12(e_2+e_3-e_5-e_7)$ & $\tfrac12(e_2-e_3+e_5+e_7)$\\
$T_4$ & $e_4$ & $\tfrac12(e_3-e_4-e_6+e_7)$ & $\tfrac12(-e_3-e_4+e_6-e_7)$\\
\bottomrule
\end{tabular}
\end{table}

Observe that the idempotents do not merely count special vertices: they select and orient an $A_2^4$ subsystem inside the Gosset polytope.  What follows requires in addition a \emph{choice}, namely that of one of the four factors as the external one; we take $T_1\cup(-T_1)$.  Its orthogonal complement in the $E_8$ root system is $E_6$.  The $240$ roots are then organized as
\begin{equation}
\boxed{240=6+72+6\cdot27,}
\label{eq:magicstarroots}
\end{equation}
which is the root-theoretic content of the Magic Star projection \cite{Truini,MCCAI23-Rosenfeld}.  At Lie-algebra level one obtains
\begin{equation}
\e_8=(\e_6,\mathbf1)\oplus(\mathbf1,\slthree)
\oplus(\mathbf{27},\mathbf3)
\oplus(\overline{\mathbf{27}},\overline{\mathbf3}).
\label{eq:magicstar}
\end{equation}
The six weights of $\mathbf3\oplus\overline{\mathbf3}$ give the six arms of the star; each arm carries a set of $27$ roots whose root spaces form a minuscule $\mathbf{27}$ or $\overline{\mathbf{27}}$ of $E_6$, the projected roots being its weight set with multiplicity one.  After choosing a Jordan basepoint, the $\mathbf{27}$ is identified with the Albert algebra $J_3(\Oa_\C):=J_3(\Oa)\otimes_\R\C$ \cite{CMZ24-MinimalCayley,CMZ25-Collineations}.  More canonically, opposite $27$-fibers form Jordan pairs, with triple product induced by the ambient Lie bracket \cite{Loos,Truini,BartonSudbery}.

The idempotent split \eqref{eq:fourtriangles} explains at the same time what remains in the central $E_6$ sector.  Once $T_1$ is used externally, the remaining nine idempotents are
\[
I_{\rm cent}=T_2\sqcup T_3\sqcup T_4.
\]
Adding their antipodes gives $18$ roots, precisely a full-rank subsystem
\[
I_{\rm cent}\cup(-I_{\rm cent})=\Phi(A_2^3)\subset\Phi(E_6).
\]
The corresponding $\Z_3$-grading is the trinification decomposition
\begin{equation}
\boxed{\e_6=\slthree^{\oplus3}\oplus(\mathbf3,\mathbf3,\mathbf3)
\oplus(\overline{\mathbf3},\overline{\mathbf3},\overline{\mathbf3}).}
\label{eq:trinification}
\end{equation}
At the group level \cite{Slansky} the connected trinification subgroup is
\begin{equation}
\frac{SU(3)\times SU(3)\times SU(3)}{\Z_3^{\rm diag}}\subset E_6.
\label{eq:globaltrin}
\end{equation}
The quotient appears because the diagonal center acts trivially on the $27$,
\begin{equation}
\mathbf{27}=(\mathbf3,\overline{\mathbf3},\mathbf1)
\oplus(\mathbf1,\mathbf3,\overline{\mathbf3})
\oplus(\overline{\mathbf3},\mathbf1,\mathbf3),
\label{eq:27trin}
\end{equation}
so that $27=9+9+9$.  In the present integral Petersson model, \eqref{eq:taucircles} cyclically permutes precisely the three idempotent triangles defining the three internal $A_2$ factors.  The nine remaining Okubo idempotents therefore provide an arithmetic skeleton for the trinification of the central $E_6$.

\begin{figure}[t]
\centering
\begin{tikzpicture}[scale=.78,transform shape,
box/.style={draw,rounded corners,align=center,inner sep=4pt,minimum width=2.35cm},
arr/.style={-{Latex[length=1.7mm]},thick}]
\node[box] (o) at (0,0) {$\mathbb O$\\$1$ idempotent\\full $E_8$};
\node[box] (p) at (4.0,0) {$p\mathbb O$\\$57=1+56$\\central $E_7$};
\node[box] (k) at (8.0,0) {Okubo\\$12=3+9$\\$A_2+E_6$};
\draw[arr] (o)--(p);
\draw[arr] (p)--(k);
\node[align=center,font=\small] at (4,-1.5) {$240$ roots fixed throughout; only the product and the idempotent polarization change.};
\end{tikzpicture}
\caption{Three multiplicative organizations of the same $E_8$ shell.}
\label{fig:three}
\end{figure}

\section{Symmetries and interpretation}

The three decompositions above involve several different symmetry groups and it is useful to keep them separate.  The real algebra automorphism groups are
\[
\Aut(\Oa)\simeq G_2,\qquad \Aut(\pO)\simeq G_2,\qquad \Aut(\Ok)\simeq PSU(3),
\]
the centre of $SU(3)$ acting trivially on $\Ok$.
These are symmetries of the multiplication on the eight-dimensional composition algebra.  By contrast, $E_7$ and $E_6$ arise as symmetries of the modules exposed by the root gradings of $\e_8$.  In the para-octonionic decomposition the $56$ carries the symplectic and quartic invariants characteristic of $E_7$.  In the Okubo Magic Star the $27$ carries the cubic invariant of $E_6$, equivalently the determinant of the complex Albert algebra.  Fixing a Jordan unit inside the $27$ reduces the structure group $E_6$ to $F_4$ \cite{SpringerVeldkamp}, whereas the nine central Okubo idempotents live in the $E_6$ transformation sector and select instead the subgroup \eqref{eq:globaltrin}.  These two kinds of idempotents should therefore not be identified.

A concise summary is given in Table~\ref{tab:summary}.  The sequence
\[
E_8\longrightarrow E_7\longrightarrow E_6
\]
should not be read as a chain produced by quotienting the same set of idempotents.  Rather, the same $E_8$ shell is equipped with three different composition products, and each idempotent equation singles out a different canonical stratification of that shell.  The octonionic unit leaves the full geometry visible; the para-idempotents identify the $E_7$ orthogonal core and its $56$-dimensional Freudenthal module; the Okubo idempotents select oriented $A_2$ data, the $E_6+A_2$ Magic Star, and then $A_2^3$ trinification inside $E_6$.

\begin{table}[t]
\centering
\caption{Exceptional structures exposed by the three integral products on $E_8$.}
\label{tab:summary}
\begin{tabular}{llll}
\toprule
Product & Idempotent pattern & Root set and module & Principal invariant\\
\midrule
$\Oa$ & $1$ & full $E_8$ shell & quadratic norm\\
$p\Oa$ & $1+56$ & $56$ roots; $\mathbf{56}$ of $E_7$ & symplectic $+$ quartic\\
Okubo & $3+9$ & $27$-root fibers; $\mathbf{27}$ of $E_6$ & cubic Jordan norm\\
\bottomrule
\end{tabular}
\end{table}

\section{Conclusions and future developments}

We have followed a single arithmetic thread from integral numbers to exceptional Lie-theoretic decompositions.  The Coxeter--Dickson order supplies the common additive $E_8$ lattice; the octonionic, para-octonionic and integral Petersson products supply three different multiplications on it.  The metric data are unchanged, but the idempotents are not.  Their counts $1$, $57$ and $12$ organize the same $240$ vertices in three distinct ways.

For the para-product the equation $x\bullet x=x$ is equivalent, on the norm-one shell, to $x^3=1$.  The $56$ nontrivial solutions lie exactly in one root layer adjacent to the distinguished root $1$, and the remaining $126$ orthogonal roots form $E_7$.  Thus the para-idempotents recover the contact decomposition of $E_8$ and expose the Freudenthal $56$.

For the Okubo product the $12$ integral idempotents form four mutually orthogonal oriented $A_2$ triangles and thus determine, with no further choice, a rank-eight subsystem $A_2^4\subset\PhiE$.  Choosing one triangle as the external $A_2$ of the Magic Star leaves nine idempotents in the central $E_6$; together with their nine antipodes they form the $18$ roots of a full-rank $A_2^3\subset E_6$, hence the trinification grading.  The six $27$-fibers of the Magic Star then connect the arithmetic picture with the Jordan pair and Albert algebra descriptions of the exceptional series.

Several questions remain open. It would be useful to characterize intrinsically the four Okubo-selected Magic Stars, to compare the para-idempotent $56$ with the integral Freudenthal geometry, and to understand how the other three Okubo-preserved lattices $C_8$, $A_2^4$ and $D_4^2$ fit into the same hierarchy of graded root decompositions.  These problems suggest that integral idempotents may provide a useful arithmetic bridge between non-unital composition algebras and exceptional Lie theory.

\subsubsection*{Acknowledgments.}
The authors thank Francesco Zucconi for helpful discussions on Okubo algebras, integral numbers, lattices and root systems.

\end{document}